\documentclass{amsart}

\usepackage{amssymb}
\usepackage{amsmath}
\usepackage[all]{xy}
\usepackage{enumerate}
\usepackage[usenames,dvipsnames]{pstricks}
\usepackage{epsfig}
\usepackage{pst-grad} 
\usepackage{pst-plot} 

\newtheorem{theorem}[subsection]{Theorem}

\newtheorem{lemma}[subsection]{Lemma}

\theoremstyle{remark}

\newtheorem{rem}[subsection]{Remark}

\theoremstyle{definition}

\DeclareMathOperator{\BOX}{Box}
\DeclareMathOperator{\Hom}{Hom}

\DeclareMathOperator{\id}{id}

\DeclareMathOperator{\orb}{orb}

\begin{document}

\title{A Geometric Interpretation of Stanley's Monotonicity Theorem}

\author{Alan Stapledon}
\address{Department of Mathematics, University of Michigan, Ann Arbor, MI 48109, U.S.A.}
\email{astapldn@umich.edu}


\thanks{The author would like to thank Professor Jiang and Professor Stanley for very useful emails.}


\begin{abstract}
We present a new geometric proof of Stanley's monotonicity theorem for lattice polytopes, using an interpretation of $\delta$-polynomials of lattice polytopes in terms of orbifold Chow rings. 
\end{abstract}

\maketitle

\section{Introduction}

Let $P$ be a $d$-dimensional lattice polytope in a lattice $N$ of rank $n$. That is, $P$ is the convex hull of finitely many points in $N \cong \mathbb{Z}^{n}$. 
If $m$ is a positive integer, then let 
$f_{P}(m) := \# \left( mP \cap N \right)$ denote the 
number of lattice points in the $m$'th dilate of $P$. A famous theorem of Ehrhart \cite{ehrhartpolynomial} asserts that
$f_{P}(m)$ is a polynomial in $m$ of degree $d$, called the 
\emph{Ehrhart polynomial} of $P$.
The generating series of $f_{P}(m)$ can be written in the form
\begin{equation*}
  \frac{ \delta_{P}(t) }{ (1 - t)^{d + 1} } = 
\sum_{m \geq 0} f_{P}(m) \, t^{m}  \, ,
\end{equation*}
where $\delta_{P}(t) = \delta_{0} + \delta_{1}t + \cdots +  \delta_{d}t^{d}$ is a polynomial of degree at most $d$ with 
integer coefficients,
called the \emph{$\delta$-polynomial} of $P$. 
Using techniques from commutative algebra, Stanley proved that 
the coefficients $\delta_{i}$ are non-negative \cite{StaDecompositions} and proved that $\delta$-polynomials of lattice polytopes satisfy the following monotonicity property \cite[Theorem 3.3]{StaMonotonicity}. 
An alternative combinatorial proof of these results was given by Beck and Sottile in \cite{BSIrrational}.
If $f(t) = \sum_{i} f_{i} t^{i}$ and $g(t) = \sum_{i} g_{i} t^{i}$ are polynomials with integer coefficients, then we write $f(t) \leq g(t)$ if $f_{i} \leq g_{i}$ for all $i \geq 0$. 

\begin{theorem}[Stanley's Monotonicity Theorem]\label{main}
If $Q \subseteq P$ are lattice polytopes, 
 then $\delta_{Q}(t) \leq \delta_{P}(t)$. 
\end{theorem}




We now present a new geometric proof of Stanley's theorem. 
We first recall the following geometric interpretation of $\delta$-polynomials of lattice polytopes.
After replacing $N$ with its intersection with the affine span of $P$, we may assume that $N$ has rank $d$.
Let $\mathcal{T}$ be a regular, lattice triangulation of $P$ and let $\sigma$ denote the cone over $P \times \{ 1\}$ in $N_{\mathbb{R}} \times \mathbb{R}$, where $N_{\mathbb{R}} = N \otimes_{\mathbb{Z}} \mathbb{R}$. The triangulation $\mathcal{T}$ induces a simplicial fan refinement $\triangle$ of $\sigma$, with cones given by the cones over the faces of $\mathcal{T}$, and we may consider  the $(d + 1)$-dimensional, simplicial toric variety $Y = Y(\triangle)$
associated to $\triangle$. The toric variety $Y$ is \emph{semi-projective} in the sense that it contains a torus-fixed point and is projective over its affinisation $Y(\sigma)$ \cite{HSToric}. 
The cohomology ring $H^{*}(X, \mathbb{Q})$ of a semi-projective, simplicial toric variety $X$ was computed by Hausel and Sturmfels in \cite{HSToric}, and it was observed by Jiang and Tseng \cite[Lemma 2.7]{JTNote} that Hausel and Sturmfel's proof, along with the results in \cite[Section 5.1]{FulIntroduction}, imply that $H^{*}(X, \mathbb{Q})$ is isomorphic to  the Chow ring $A^{*}(X, \mathbb{Q})$. 

The orbifold Chow ring of a Deligne-Mumford stack was introduced by Abramovich, Graber and Vistoli \cite{AGVAlgebraic} as the algebraic analogue of Chen and Ruan's orbifold cohomology ring \cite{CRNew}. Borisov, Chen and Smith introduced the notion of a toric stack in \cite{BCSOrbifold} and showed that any simplicial, semi-projective toric variety $X$ has the canonical structure of a Deligne-Mumford stack. The orbifold Chow ring $A^{*}_{\orb}(X, \mathbb{Q})$ of $X$ is a $\mathbb{Q}$-graded $\mathbb{Q}$-algebra and was computed by 
Jiang and Tseng in \cite{JTNote}, generalising results in \cite{BCSOrbifold} (Remark \ref{explicit}).
The following combinatorial observation follows from \cite[Theorem 4.6]{YoWeightI} (c.f. \cite[Corollary 1.2]{KarEhrhart}).
\begin{equation}\label{express}
\delta_{P}(t) = \sum_{i \in \mathbb{Q}} \dim_{\mathbb{Q}} A^{i}_{\orb}(Y, \mathbb{Q})t^{i}.
\end{equation}



If $Q$ is a lattice polytope contained in $P$, then let $N'$ denote the intersection of $N$ with the affine span of $Q$ and let $\sigma'$ denote the cone over $Q \times \{ 1 \}$ in $(N')_{\mathbb{R}} \times \mathbb{R}$. One verifies that we may choose a regular, lattice triangulation $\mathcal{T}$ of $P$ which restricts to a regular, lattice triangulation of $Q$. In this case, 
the fan $\triangle$ refining $\sigma$ restricts to a fan $\Sigma$ refining $\sigma'$ and we may consider the semi-projective 
toric variety
$Y' = Y'(\Sigma)$. The inclusion of $N'$ in $N$ induces a locally closed toric immersion  $Y' \hookrightarrow Y$ and
a restriction map between the corresponding orbifold Chow rings.  We will prove the following lemma in Section \ref{details}.  

\begin{lemma}\label{open}
The morphism $Y' \hookrightarrow Y$ induces a surjective graded ring homomorphism $A^{*}_{\orb}(Y, \mathbb{Q}) \rightarrow A^{*}_{\orb}(Y', \mathbb{Q})$. 
\end{lemma}

By (\ref{express}),
$\delta_{P}(t) = \sum_{i \in \mathbb{Q}} \dim_{\mathbb{Q}} A^{i}_{\orb}(Y, \mathbb{Q})t^{i}$ and
$\delta_{Q}(t) = \sum_{i \in \mathbb{Q}} \dim_{\mathbb{Q}} A^{i}_{\orb}(Y', \mathbb{Q})t^{i}$, and we conclude that $\delta_{Q}(t) \leq \delta_{P}(t)$, as desired. 



\begin{rem}
If we regard the empty face as a face of the triangulation $\mathcal{T}$ of dimension $-1$, then the $h$-vector of $\mathcal{T}$ is defined by
\begin{equation*} 
h_{\mathcal{T}}(t) = \sum_{F} t^{\dim F + 1} (1 - t)^{d - \dim F},
\end{equation*}
where the sum ranges over all faces $F$ in $\mathcal{T}$. It is a well known fact that $0 \leq h_{\mathcal{T}}(t) \leq \delta_{P}(t)$ and $h_{\mathcal{T}}(t) = \delta_{P}(t)$ if and only if $\mathcal{T}$ is a unimodular triangulation \cite{BMLattice, PayEhrhart}. We have the following geometric interpretation of this result. 

It follows from the definition of the orbifold Chow ring (see Section 2) that $A^{*}(Y, \mathbb{Q})$ is a direct summand of $A^{*}_{\orb}(Y, \mathbb{Q})$ and $A^{*}(Y, \mathbb{Q}) = A^{*}_{\orb}(Y, \mathbb{Q})$ if and only if $Y$ is smooth. The result now follows from the fact that 
$h_{\mathcal{T}}(t) = \sum_{i \geq 0} \dim_{\mathbb{Q}} A^{i}(Y, \mathbb{Q}) t^{i}$ \cite[Corollary 2.12]{HSToric} and the fact that $Y$ is smooth if and only if $\mathcal{T}$ is a unimodular triangulation. 

\end{rem}




All varieties and stacks will be over the complex numbers. In Section \ref{details}, we will review orbifold Chow rings and prove Lemma \ref{open}. 

\section{Orbifold Chow Rings}\label{details}

The orbifold Chow ring $A^{*}_{\orb}(\mathcal{X}, \mathbb{Q})$ of a Deligne-Mumford stack $\mathcal{X}$ was introduced by Abramovich, Graber and Vistoli as the degree $0$ piece of the small quantum cohomology ring of $\mathcal{X}$ \cite{AGVAlgebraic}. We will review the structure of 
$A^{*}_{\orb}(\mathcal{X}, \mathbb{Q})$ as a $\mathbb{Q}$-graded vector space and refer the reader to \cite{AGVAlgebraic} for the relevant details and the description of the ring structure of $A^{*}_{\orb}(\mathcal{X}, \mathbb{Q})$. The \emph{inertia stack} $\mathcal{I}\mathcal{X}$ of $\mathcal{X}$ is a Deligne-Mumford stack whose objects consist of pairs $(x, \alpha)$, where $x$ is an object of $\mathcal{X}$ and $\alpha$ is an automorphism of $x$. If $\mathcal{X}_{1}, \ldots, \mathcal{X}_{r}$ denote the connected components of $\mathcal{I}\mathcal{X}$, then
\begin{equation*}
A_{\orb}^{*}(\mathcal{X}, \mathbb{Q}) = \bigoplus_{j = 1}^{r} A^{*}(|\mathcal{X}_{j}|, \mathbb{Q})[s_{j}],
\end{equation*}
where $|\mathcal{X}_{j}|$ is the coarse moduli space of $\mathcal{X}_{j}$,
$s_{j} \in \mathbb{Q}$ is the \emph{age} of $\mathcal{X}_{j}$ and $[s_{j}]$ denotes a grading shift by $s_{j}$. If we identify $\mathcal{X}$ as the connected component of $\mathcal{I} \mathcal{X}$ whose objects consist of pairs $(x, \id)$, where $x$ is an object of $\mathcal{X}$ and $\id$ is the identity automorphism of $x$, then the age of  $\mathcal{X}$ is $0$ and $A^{*}(|\mathcal{X}|, \mathbb{Q})$ is a direct summand of $A^{*}_{\orb}(\mathcal{X}, \mathbb{Q})$. 

Continuing with the notation of the introduction, recall that $P$ is a $d$-dimensional lattice polytope in  a lattice $N$ of rank $d$ and $\mathcal{T}$ is a regular lattice triangulation of $P$. Recall that $\mathcal{T}$ induces a fan refinement $\triangle$ of the cone $\sigma$ over $P \times \{1\}$ in $N_{\mathbb{R}} \times \mathbb{R}$, and that $Y = Y(\triangle)$ is the associated $(d + 1)$-dimensional, semi-projective, simplicial toric variety. 
There is a canonical Deligne-Mumford stack $\mathcal{Y}$ with coarse moduli space $Y$ \cite{BCSOrbifold}. 
If $F$ is a non-empty face of $\mathcal{T}$ with vertices $v_{1}, \ldots, v_{s}$, then set
\begin{equation*}
\BOX(F) = \{ w \in N_{\mathbb{R}} \times \mathbb{R} \mid  w = \sum_{i = 1}^{s} q_{i}(v_{i}, 1) \textrm{  for some  }
0 < q_{i} < 1 \},
\end{equation*}
and let  $\BOX(\emptyset) = \{ 0  \in N_{\mathbb{R}} \times \mathbb{R} \}$. Borisov, Chen and Smith \cite{BCSOrbifold} showed that the 
inertia stack of $\mathcal{Y}$ decomposes into connected components as 
$\mathcal{I}\mathcal{Y} = \coprod_{F \in \mathcal{T}} \coprod_{w \in BOX(F) \cap (N \times \mathbb{Z})} \mathcal{Y}_{w}$, where $\mathcal{Y}_{w} = \mathcal{Y}$ if $w = 0$ and, if $w \neq 0$, then $|\mathcal{Y}_{w}|$ is isomorphic to the torus-invariant subvariety $V(F)$ of $Y$ corresponding to the cone over $F \times \{1\}$ in $\triangle$.
Moreover, if $\psi : N_{\mathbb{R}} \times \mathbb{R} \rightarrow \mathbb{R}$ denotes projection onto the second co-ordinate, then 
the age of $\mathcal{Y}_{w}$  is  $\psi(w) \in \mathbb{Z}$.

Recall that if $Q$ is a lattice polytope contained in $P$, then $N'$ is the intersection of $N$ with the affine span of $Q$ and 
the fan $\triangle$ restricts to a fan $\Sigma$ refining the cone $\sigma'$ over $Q \times \{ 1 \}$ in $(N')_{\mathbb{R}} \times \mathbb{R}$. 
If $\mathcal{Y}'$ denotes the canonical Deligne-Mumford stack with coarse moduli space $Y' = Y'(\Sigma)$, then 
the inclusion of $N'$ in $N$ induces an inclusion of $\mathcal{Y}'$ as a closed substack of $\mathcal{Y}' \times (\mathbb{C}^{*})^{\dim P - \dim Q}$ and an inclusion of $\mathcal{Y}' \times (\mathbb{C}^{*})^{\dim P - \dim Q}$ as an open substack of $\mathcal{Y}$. These inclusions induce a corresponding 
restriction map $\iota: A^{*}_{\orb}(\mathcal{Y}, \mathbb{Q}) \rightarrow A^{*}_{\orb}(\mathcal{Y}', \mathbb{Q})$, which we describe below (c.f. Remark \ref{explicit}). 

If $\mathcal{T}|_{Q}$ denotes the restriction of $\mathcal{T}$ to $Q$, then the inertia stack of $\mathcal{Y}'$ decomposes into connected components as 
$\mathcal{I}\mathcal{Y}' = \coprod_{F \in \mathcal{T}|_{Q}} \coprod_{w \in BOX(F) \cap (N \times \mathbb{Z})} \mathcal{Y}'_{w}$, where $\mathcal{Y}'_{w} = \mathcal{Y}'$ if $w = 0$ and, if $w \neq 0$, then 
the age of $\mathcal{Y}'_{w}$  is  $\psi(w)$ and 
$|\mathcal{Y}'_{w}|$ is isomorphic to the torus-invariant subvariety $V(F)'$ of $Y'$ corresponding to the cone over $F \times \{1\}$ in $\Sigma$. For each face $F \in \mathcal{T}|_{Q}$, 
the inclusion of $N'$ in $N$ induces a closed immersion  $V(F)' \hookrightarrow V(F)' \times (\mathbb{C}^{*})^{\dim P - \dim Q}$ and 
an open immersion  $V(F)' \times (\mathbb{C}^{*})^{\dim P - \dim Q} \hookrightarrow V(F)$. 
The corresponding restriction map $\nu_{F}: A^{*}(V(F), \mathbb{Q}) \rightarrow A^{*}(V(F)', \mathbb{Q})$ is surjective since if 
$W'$ is an irreducible closed subvariety of $V(F)'$ and $W$ denotes the closure of $W' \times (\mathbb{C}^{*})^{\dim P - \dim Q}$ in $V(F)$, then 
$\nu_{F}([W]) = [W']$. The restriction map $\iota: A^{*}_{\orb}(\mathcal{Y}, \mathbb{Q}) \rightarrow A^{*}_{\orb}(\mathcal{Y}', \mathbb{Q})$ has the form 
\[
\iota: \coprod_{F \in \mathcal{T}} \coprod_{w \in BOX(F) \cap (N \times \mathbb{Z})} A^{*}(|\mathcal{Y}_{w}|, \mathbb{Q})[\psi(w)] \rightarrow 
\coprod_{F \in \mathcal{T}|_{Q}} \coprod_{w \in BOX(F) \cap (N \times \mathbb{Z})} A^{*}(|\mathcal{Y}'_{w}|, \mathbb{Q})[\psi(w)],
\]
where for each $F \in \mathcal{T}$ and $w \in BOX(F) \cap (N \times \mathbb{Z})$, $\iota$ restricts to $\nu_{F}$ (with a grading shift) on $A^{*}(|\mathcal{Y}_{w}|, \mathbb{Q})[\psi(w)]$ if $F \subseteq Q$ and restricts to zero otherwise. 
One can verify from the description of the ring structure of an orbifold Chow ring in \cite{AGVAlgebraic} that $\iota$ is a ring homomorphism. 
We conclude that 
$\iota$ is a surjective ring homomorphism, thus establishing  Lemma \ref{open}.

\begin{rem}\label{analogue}
The dimensions of the graded pieces of $A^{*}(V(F), \mathbb{Q})$ are equal to the coefficients of an $h$-vector  of a fan  \cite[Corollary 2.12]{HSToric}. 
The analogous combinatorial proof of Stanley's theorem goes as follows: 
one can express $\delta_{P}(t)$ and $\delta_{Q}(t)$ as sums of shifted $h$-vectors  \cite{BMLattice, PayEhrhart}, and then apply Stanley's monotonicity theorem for $h$-vectors \cite{StaMonotonicity} to conclude the result. 
\end{rem}

\begin{rem}\label{explicit}
Consider the deformed group ring $\mathbb{Q}[N \times \mathbb{Z}]^{\triangle} := \oplus_{v \in \sigma \cap (N \times \mathbb{Z})} \mathbb{Q} \cdot y^{v}$, with ring structure defined by
\begin{displaymath}
y^{v} \cdot y^{w} =  \left\{ \begin{array}{ll}
y^{v + w} & \textrm{if there exists a cone  } \tau \in \triangle \textrm{containing} v \textrm{and} w\\
0 & \textrm{otherwise}.
\end{array} \right.
\end{displaymath}
If $v_{1}, \ldots, v_{t}$ denote the vertices of $\mathcal{T}$ and $M = \Hom_{\mathbb{Z}}(N, \mathbb{Z})$, then Jiang and Tseng \cite[Theorem 1.1]{JTNote}
showed that there is an isomorphism of rings
\[
A^{*}_{\orb}(Y, \mathbb{Q}) \cong \frac{\mathbb{Q}[N \times \mathbb{Z}]^{\triangle}}{ \{ \sum_{i = 1}^{t} \langle (v_{i}, 1), u  \rangle y^{(v_{i}, 1)} \mid u \in M \times \mathbb{Z} \}}. 
\]
Similarly, if $v_{1}, \ldots, v_{s}$ are the vertices of $\mathcal{T}|_{Q}$ and $M' =  \Hom_{\mathbb{Z}}(N', \mathbb{Z})$, then
\[
A^{*}_{\orb}(Y', \mathbb{Q}) \cong \frac{\mathbb{Q}[N' \times \mathbb{Z}]^{\Sigma}}{ \{ \sum_{i = 1}^{s} \langle (v_{i}, 1), u  \rangle y^{(v_{i}, 1)} \mid u \in M' \times \mathbb{Z} \}}. 
\]
Consider the surjective ring homomorphism $j: \mathbb{Q}[N \times \mathbb{Z}]^{\triangle} \rightarrow \mathbb{Q}[N' \times \mathbb{Z}]^{\Sigma}$
satisfying $j(y^{v}) = y^{v}$ if $v \in \Sigma$ and $j(y^{v}) = 0$ if $v \notin \Sigma$. The induced ring homomorphism 
\[
\frac{\mathbb{Q}[N \times \mathbb{Z}]^{\triangle}}{ \{ \sum_{i = 1}^{t} \langle (v_{i}, 1), u  \rangle y^{(v_{i}, 1)} 
\mid u \in M \times \mathbb{Z} \} }
\longrightarrow \frac{\mathbb{Q}[N' \times \mathbb{Z}]^{\Sigma}}{ \{ \sum_{i = 1}^{s} \langle (v_{i}, 1), u  \rangle y^{(v_{i}, 1)} 
\mid u \in M' \times \mathbb{Z} \}}
\]
corresponds to the restriction map $\iota: A^{*}_{\orb}(Y, \mathbb{Q}) \rightarrow A^{*}_{\orb}(Y', \mathbb{Q})$ under the above isomorphisms. The existence of such a ring homomorphism was used by Stanley in his original commutative algebra proof of Theorem \ref{main} \cite{StaMonotonicity}. 
\end{rem}




\bibliographystyle{amsplain}
\bibliography{alan}

\end{document}